
\magnification=\magstep1
\documentstyle{amsppt}

\def\pr{\noindent{\sl Proof. }}
\def\re{\noindent{\sl Remark. }}
\def\E{{\Cal E}}

\def\H{{\Bbb H}}
\def\I{{\Cal I}}

\def\K{{\Cal K}}
\def\L{{\Cal L}} 
\def\M{{\Cal M}} 

\def\O{{\Cal O}}
\def\P{{\Bbb P}}
\def\Z{{\Bbb Z}}
\def\Tors{{\operatorname {Tors}}}
\def\Ext{{\operatorname{Ext}}}
\def\coker{{\operatorname{coker}}}

\def\hom{{\Cal H}om}
\def\ra{\longrightarrow}
\def\raa{\rightarrow}

\def\ot{\otimes}
\def\op{\oplus}

\def\today{\ifcase\month \or January \or February \or March \or April
\or May \or June \or July \or August \or September \or October
\or November \or December \fi \space \number\day,  \number \year}
\def\Ann{\operatorname{{Ann}}}
\def\Hom{{\operatorname{Hom}}}
\def\:{\colon}

\topmatter
\title 
Pluricanonical systems on surfaces with small $K^2$.
\endtitle
\author Adrian Langer
\endauthor

\rightheadtext{pluricanonical systems}
{\date{\today}\enddate}

\address{ Adrian Langer:
Instytut Matematyki UW, 
ul.~Banacha 2, 02--097 Warszawa, Poland}
\endaddress
\email{alan\@mimuw.edu.pl}\endemail

\abstract{}
\endabstract

\thanks The author was partially supported by Foundation for Polish Science
\endthanks


\endtopmatter

\document 
\vskip .3 cm

\heading  Introduction\endheading

The main purpose of this paper is studying of the properties of pluricanonical
systems on surfaces with $K^2\le 4$.

The first part is devoted to the study of the spannedness of 
bicanonical systems on surfaces of general type with small $K^2$
or, which is the same, on canonical models of such surfaces.

The study of the  bicanonical system on a canonical surface has two
advantages. Namely, it is well known that linear systems
tend to have less singular points as base points (e.g., see Theorem 1.1).
The second advantage is that studying of an ample line bundle is easier 
than that of a nef one.

The main result of this part is the following

\proclaim{Theorem 0.1}
Let $X$ be a canonical surface with $K_X^2=4$.
Then $|2K_X|$ has no fixed part and every base point of $X$ is smooth.
In particular, the bicanonical system on the minimal model of $X$ has no base
component.
\endproclaim

A well-known theorem says that $|2K_X|$ is base point free if $K_X^2>4$
or $p_g(X)>0$ but our knowledge of base points, or even
base components, in the  case when $2\le K_X^2\le 4$
and $p_g=0$ is very limited (see e.g.~[CT] and [We]).

To prove Theorem 0.1 we use our earlier work on adjoint linear series
([La1] and [La2]) generalizing Reider's criterion to normal surfaces
to exclude singular base points and all but one possibility for the fixed
component and then we use the bilinear map lemma and 
Clifford's theorem for arbitrary Cohen--Macaulay curves (see Section 2)
to exclude the remaining case. 

\medskip

In the second part of the paper we describe which degree $2$ clusters are
contracted by $\varphi_{4K_X}$ on a numerical Godeaux surface and by
$\varphi_{3K_X}$ on a Campedelli surface.  Here we have the following
theorems (looking at them one should remember that the torsion group
$\Tors X$ is finite; see 1.6):

\proclaim{Theorem 0.2}
Let $X$ be a canonical numerical Godeaux surface. 
If $|4K_X|$ contracts a degree $2$ cluster $\zeta$, then
$\zeta$ is contained in a curve $D\in |K_X+\tau|$, where $0\ne \tau\in \Tors X$
and the morphism $\varphi _{4K_X} $ restricted to $D$ is either:
\item{1.} an embedding apart from the contracted cluster $\zeta$
if $2\tau\ne 0$ in $\Tors X$,or
\item{2.} a double covering of $\P ^1$, if $2\tau =0$ in $\Tors X$.
\endproclaim

\proclaim{Theorem 0.3}
Let $X$ be a canonical model of a numerical Campedelli surface.
If a degree $2$ cluster $\zeta$ is contracted by 
$|3K_X|$, then one of the following holds:
\item{1.}
There exists an honestly hyperelliptic curve $C\in |2K_X|$ containing
$\zeta$ and $\varphi_{3K_X}|_C$ is a double covering of $\P^1$.
\item{2.}
$\zeta$ is a scheme-theoretic intersection of two curves
$B_1\in |K_X+\tau|$ and $B_2\in |K_X-\tau|$, for some
$\tau\in \Tors X$ such that $2\tau\ne 0$.
\endproclaim

After writing down this paper (apart from the part after Lemma 5.4,
which was inspired by the paper of Kotschick)
we learnt about a related paper of Kotschick [Ko], 
who however considered mainly the case of torsion free numerically
Godeaux or Campedelli surfaces. Our Theorems 0.2 and 0.3 (together with
Proposition 5.5 giving a topological interpretation of curves appearing in
case 1 of Theorem  0.3) generalize Theorems 1 and 2, [ibid]. 

He also proved  ([Theorem 3, ibid]) a criterion for spannedness of $|2K_X|$ 
in terms  of the fundamental group of $X$
and the second Stifel-Whitney class of the tangent bundle of $X$. 
Since even the fundamental group is known only in explicit examples 
and can vary (and for some examples it is not finite), 
the criterion seems to be only of theoretical interest.

\heading 1. Preliminaries \endheading

All varieties are assumed to be defined over ${\Bbb C}$ 
(except in Section 2). In this section we will state some results 
generalizing Reider's criterion. First let us recall a special case of
Corollary 5.1.4, [La2]:

\proclaim{Theorem 1.1} ([La2], Corollary 5.1.4).
Let $X$ be a normal projective surface with only  quotient singularities
and $L$ be a nef Weil divisor on $X$ such that $K_X+L$ is Cartier. 
Assume that $L^2>{4\over r}$, where $r$ is
the order of the local fundamental group around a given point $x$.
Then $K_X+L$ is not globally generated at $x$ if and only if there exists
a connected curve $D$ containing $x$ such that $\O _D(K_X+L)$ is not
globally generated at $x$. Moreover, one can choose $D$ satisfying
the following conditions:
\item{1.} There exists an injection 
$m_x\O _D(K_X+L)\hookrightarrow \omega _D$,
\item{2.} $LD-{1\over r}\le D^2\le {(LD)^2\over L^2}$ and $0\le LD<{2\over r}$,
\item{3.} If $X$ has only Du Val singularities then 
$m_x\O _D(K_X+L)\simeq \omega _D$. In particular, $2p_aD=D(K_X+L)+1$.
\endproclaim

{\it Sketch of the proof.} 
If $K_X+L$ is not globally generated at $x$ then by Serre's construction
there exists a rank $2$ reflexive sheaf $\E\in \Ext ^1(m_x\O(K_X+L),\omega_X)$.
This sheaf is locally at the point $x$ isomorphic with $\Omega_X^{**}$.
Using this fact the theorem can be proved by using Bogomolov's
instability theorem for surfaces with only quotient singularities
(see [Ka, Lemma 2.5]).

The details of the proof together with a generalization of the theorem
can be found in [La2] or in a forthcoming paper of the author.

\proclaim{Corollary 1.2}
Let $x$ be a singular point of the canonical surface $X$. Then 
$|2K_X|$ is globally generated at $x$ unless
\item{1.} $K_X^2=1$ and $x$ is of type $A_1$, $A_2$, $A_3$, or
\item{2.} $K_X^2=2$ and $x$ is of type $A_1$.
\endproclaim

A special case of this corollary was proved in [We].

\proclaim{Definition 1.3}
Assume that a normal surface $X$ has rational singularity at $x$ and let
$f\:Y\to X$ be the minimal resolution at $x$. Let $Z$
denote the fundamental cycle (or the exceptional divisor if $x$ is smooth).
The Seshadri constant of a nef divisor $L$ at $x$ 
is the real number
$$\epsilon (L,x)=\sup\{\epsilon\ge 0| f^*L-\epsilon\cdot Z{\hbox{ is nef}}\}.$$
\endproclaim

In the proof of the next proposition we will need the following

\proclaim{Lemma 1.4}
Let $D$ be a Weil divisor on a normal surface. If $D^2\ge 0$ and $DL>0$
for some nef divisor $L$, then $D$ is pseudoeffective. If moreover
$D$ is not big then $D^2=0$ and $D$ is nef.
\endproclaim

\pr 
First part of the lemma follows easily from the Hodge index theorem.
The second one follows from the Zariski decomposition for $D$, Q.E.D.

\proclaim{Proposition 1.5} 
Let $X$ be a canonical surface and $x$ a base point of $|2K_X|$. 
Then
\item{1.} If $K_X^2=4$, then $x$ is smooth and  $\epsilon (K_X,x)=2.$
\item{2.} If $K_X^2=2$ and $x$ is singular, then $x$ is of type $A_1$
and $\epsilon (K_X,x)=1$.
\endproclaim

\pr This is just a simple corollary to Corollary 1.2 and 
[La3, Corollary 3.2].

\medskip
\noindent
{\bf 1.6.} Let $X$ be a canonical surface with $p_g=q=0$.
Recall that $X$ is called numerical Godeaux (numerical Campedelli) if 
$K_X^2=1$ ($K_X^2=2$, respectively).  Note that usually these notions are
defined for a minimal model but as we said before it is more convenient
to use the canonical model.

The torsion group
$\Tors X$ (i.e., a torsion subgroup of $\hbox{Pic} X$; it is the same for $X$
and its minimal model) of a numerical Godeaux or Campedelli surface
is finite. In fact if $X$ is numerical Godeaux then $\Tors X$
is cyclic of order at most $5$ (and all these possibilities occur).
If $X$ is numerical Campedelli then instead one can consider $\pi_1^{alg}(X)$
which is a related group. This group is finite of order at most $9$.
Moreover, if $|\pi_1^{alg}(X)|=9$ then 
$\Tors X=\pi_1^{alg}(X)=\pi_1(X)=\Z_3\op \Z_3$. If $|\pi_1^{alg}(X)|=8$
then $\pi_1^{alg}(X)=\pi_1(X)$ and we have the following possibilities:
\item{a)} $\pi_1(X)=\Tors X=\Z_2\op\Z_2\op\Z_2$, $\Z_2\op\Z_4$ or $Z_8$,
\item{b)} $\pi_1(X)=\H$ (the quaternion group) and $\Tors X=\Z_2\op\Z_2$,

\noindent
and all of them occur (this follows from explicit description of those
surfaces in [Re]).

\medskip
\noindent
{\bf 1.7.} Let us recall (see e.g., [CFHR]) that a Gorenstein curve $C$
is called {\sl honestly hyperelliptic} if there exists a finite morphism
$\varphi \: C\to \P^1$ of degree $2$. Clearly every irreducible reduced
curve of genus $2$ is honestly hyperellipic.

\heading 2. Clifford's lemma \endheading

Let us recall the following well-known bilinear map lemma:

\proclaim{Lemma 2.1} (Hopf, [Ha, Lemma 5.1])
Let $\phi :V_1\times V_2\ra W$ be a bilinear map of nonzero finite-dimensional
vector spaces (over an algebraically closed field $k$), which is 
nondegenerate, i.e., for each $v_1\ne 0$ in $V_1$ and
each $v_2\ne 0$ in $V_2$, $\phi (v_1,v_2)\ne 0$. Then
$$dim\, W\ge dim\, V_1+dim\, V_2-1.$$
\endproclaim

\proclaim{Corollary 2.2}
Let $X$ be an integral scheme defined over an algebraically closed field $k$.
Let $\L$ and $\M$ be two coherent subsheaves of the sheaf of total quotient
rings $\K _X$ such that $h^0(\L)\ge 1$ and $h^0(\M )\ge 1$.
Then
$$h^0(\L\cdot \M)\ge h^0(\L)+h^0(\M )-1,$$
where $\L \cdot \M$ is the product of sheaves in $\K _X$.
\endproclaim

\pr
We have a natural multiplication map $\L (X)\times \M (X)\ra 
(\L \cdot \M ) (X)$. If $s_1\in \L (X)$ and $s_2\in \M (X)$ then 
$s_1 s_2 \ne 0$ in $(\L \cdot \M)(X)\subset \K (X)$ and we can apply the
bilinear map lemma, Q.E.D. 
\medskip

\proclaim{Definition 2.3} A sheaf $L$ on a scheme $X$ is called invertible in 
codimension $0$, or generically invertible, if it is locally isomorphic 
to $\O _X$ at every generic point of $X$.
\endproclaim

\proclaim{Lemma 2.4} (Clifford).
Let $C$ be a Cohen-Macaulay, projective curve over an algebraically 
closed field $k$. Let $\L$ be a generically invertible 
torsion free coherent sheaf
such that $\deg \L\le 2p_aC+h^0(\O_C)-1$. 
Then there exists a generically Gorenstein subcurve 
$B\subset C$ such that $h^0(B,\L|_B)=0$ or
$$2h^0(B,\L |_B)\le \deg \L|_B+h^0(\O _B)+1.$$
\endproclaim

\pr
By the definition of degree and Serre duality we have:
$$\chi (\L)=h^0(\L) -h^0(\hom (\L , \omega _C))= \deg \L+\chi (\O _C)$$
Now if $h^0(\hom (\L , \omega _C))=0$, then
$2h^0(\L)=2\deg \L +2\chi (\O_C)\le \deg \L+h^0(\O _C)+1$,
because $\deg \L\le 2p_aC+h^0(\O_C)-1$.

Hence we can assume that $h^0(\hom (\L , \omega _C))\ne 0$ and 
$h^0(\L)\ne 0$. There exists a subcurve $B\subset C$ such that 
$\Hom (\L, \omega _C)\ne 0$ and every non-zero homomorphism 
$\varphi \:\L|_B \to \omega _B$ is generically
onto. Indeed, if $\psi$$\in \Hom (\L, \omega _C)$ is not generically onto
then we can choose a subcurve $B'\subset C$ defined by $\Ann \psi$, such
that $\psi$ has factorisation 
$$\L\to \L|_{B'}\mathop{\to}^{\psi_{B'}} \omega_{B'}\subset \omega_C,$$  
where $\psi_{B'}$ is generically onto (see [CFHR, Lemma 2.4]).
If necessary we continue this process for other homomorphisms
$\L|_{B'}\to \omega_{B'}$ until we get the required curve. 

Obviously we can assume that $h^0(B,\L|_B)\ne 0$. Then we have a natural 
pairing  $H^0(B,\L |_B)\times \Hom (\L|_B,\omega _B)\to
H^0(\omega_B)$ and by the above we see that assumptions of the
bilinear map lemma are satisfied. Therefore
$$h^0(B,\L|_B)+h^0(\hom (\L|_B , \omega _B))\le 
h^0(\omega _B)+1.$$
Hence 
$$2h^0(B,\L|_B) -\deg \L|_B-\chi (\O_B)=h^0(B,\L|_B)+h^0(\hom (\L|_B , 
\omega _B))\le h^0(\omega _B)+1,$$
which proves the lemma.
\medskip

We will usually apply Clifford's lemma for rank one torsion-free
sheaves on a reduced, irreducible curve.

\heading 3. Proof of Theorem 0.1 \endheading

The second part of the theorem is contained in Corollary 1.2.

Let us write $|2K_X|=|M| +V$, where $V$ is a fixed part.
We can assume that $p_g(X)=q(X)=0$, since otherwise $2K_X$ is
base point free.
Then $|M|$ is not composed with a pencil (see [Xi]).
Therefore a general member of $|M|$ is irreducible and reduced.
By an abuse of notation we will write $M$ for
the general member of $|M|$. 

\proclaim{Lemma 3.1}
$$2K_X^2\le p_aM\le {1\over 2}(K_XM+M^2)+1$$
\endproclaim

\pr
From the long cohomology exact sequence corresponding to
$$0\ra \O_X(K_X-M) \ra \O (K_X)\ra \O _M(K_X)\ra 0,$$
one can get 
$$0\ra H^1(K_X|_M)\ra H^2(K_X-M)\ra H^2(K_X)\ra 0.$$
Hence by Serre duality we have
$$h^1(K_X|_M)= h^2(K_X-M)- h^2(K_X)=h^0(M)-1=h^0(2K_X)-1=K_X^2.$$
Similarly one can compute 
$$h^1(2K_X|_M)=h^2(2K_X-M)- h^2(2K_X)=h^0(K_X-V)=0$$ 
since $h^0(K_X)=0$ by assumption.
Now using the Riemann--Roch theorem on $M$ we get:
$$h^0(K_X|_M)=\chi (K_X|_M)+h^1(K_X|_M)=\chi (\O_M)+K_XM+K_X^2$$
and
$$h^0(2K_X|_M)=\chi (2K_X|_M)+h^1(2K_X|_M)=\chi (\O_M)+2K_XM.$$
Using Corollary 2.2, we obtain:
$$2h^0(K_X|_M)-1\le h^0(2K_X|_M).$$ 
Substituting the above equalities for both sides of the previous
inequality, we get 
$2K_X^2\le p_aM$. The inequality   $2p_a(M)-2\le (K_X+M)M$ 
(a ``numerical subadjunction'') is a consequence of the Riemann--Roch 
theorem for surfaces  with at  most Du Val  singularities
(see e.g., [La1, Theorem 2.1]).
In fact, from this  theorem it follows that an equality holds if and  only
if $M$ is Cartier, Q.E.D.
\medskip

{\re }
The lemma works also for surfaces $K_X^2=2$ or $3$ and $p_g(X)=0$.
It limits the number and intersection numbers of the possible base components
of the bicanonical system.

\medskip

From Lemma 1.4 it follows that $K_XV\ge 2\, mult_xV$ for any $x\in V$.
In particular $K_XV\ge 2$. 

Since $V$ does not pass through the singular points  of $X$ by Corollary 1.2,
$M=2K_X-V$ is a Cartier  divisor and hence $2p_aM-2=(K_X+M)M$. 
Using this equality and Lemma 3.1  we get  $M^2\ge 14-KM=6+KV\ge 8$.
By the Hodge index theorem $(KM)^2\ge K^2M^2\ge 4\cdot  8$, hence $KM\ge 6$.
Using $KM+KV=2K^2=8$, we are left with only one case:
$KV=2$,  $KM=6$. Since $M^2=2p_aM-2-K_XM$ is divisible by $2$  and
$8\le M^2\le (KM)^2/K^2=9$, we get $M^2=8$, $MV=4$ and $V^2=0$.

Now it is easy to see that there is only
one possibility: $K_XV=2$, $K_XM=6$, $M^2=8$, $MV=4$, $V^2=0$.

\medskip
In this case we have a sequence:
$$0\ra \O (M-V)\ra \O (M)\ra \O_V(M)\ra 0 $$
and
$$h^0(\O_V (M))\le {1\over 2}MV+1=3$$
by Clifford's theorem, since $MV\le 2p_a(V)=4$. Hence we get
$h^0(\O (M-V))\ge h^0(\O (M))-h^0(\O _V(M))\ge 2.$
Let $G$ be a generic member of $|M-V|$. We will prove that 
$h^0(\O_A)\le 2$ for every subcurve $A\subset G$. 
Then, since $4=MG\le 2p_aG=6$, using Clifford's theorem, we get 
$$4=h^0(M)-h^0(V)\le h^0(\O _{A}(M))\le {1\over 2}(MA+h^0(\O_A)+1)
< {1\over 2}M(M-V)+2=4,$$ 
hence a contradiction.

Because $M-V=2(K_X-V)$, $K_X-V$ is pseudoeffective. If it is big, then
we can apply [La1, Theorem 3.6], since  $h^1(K_X+(K_X-V))=h^1(M)=
h^0(\O_V(2K_X))\ne 0$. Therefore there exist divisors $A$ and $B$ such 
that $K_X-V=A+B$, $A-B$ is pseudoeffective (and numerically nontrivial) 
and $B$ is effective. Hence $1\le BK_X<AK_X$ and we get a contradiction
with $(A+B)K_X=(K_X-V)K_X=2$. 

Hence $M-V$ is not big and it is easy to see that it is nef.
Therefore
if we write $|M-V|=|D|+F$, where $F$ is the fixed part of $|M-V|$, 
then $D(D+F)=F(D+F)=(M-V)^2=0$.
Hence $DF\le D^2+DF=0$ and we have a contradiction unless $F=0$.
In the latter case $|D|$ is composed with a pencil.
Because $q(X)=0$, $|D|$ is composed with a rational pencil. Moreover,
$|D|$ has no base points since $D^2=0$. Therefore $D=f^*\O _{{\P}^1}(2)$,
where $f\:X\to {\P}^1$ is the morphism defined by $|D|$,
and $h^0(\O _D)=2$, Q.E.D.

\heading 4. Numerical Godeaux surfaces \endheading

\proclaim{Proposition 4.1}
Let $X$ be a canonical numerical Godeaux surface.
Then $|4K_X|$ contracts a degree $2$ cluster $\zeta$ 
if and only if $\zeta$ is a scheme-theoretical intersection of two curves
$D_1\in |K_X+\tau|$ and $D_2\in |2K_X-\tau|$, where $0\ne\tau\in \Tors X$.
\endproclaim

\pr
By [La1, Theorem 0.2] applied to $L=3K_X$ one can see that 
if $\zeta$ is contracted by $|4K_X|$, then there exists an
effective Cartier divisor $D$
passing through $\zeta$ such that $K_XD=1$, $p_aD=2$, $D\equiv K_X$
and $\omega _D\simeq\I_\zeta\O_D(4K_X)$.
Moreover, from the construction of the divisor $D$, we see that
the  bundle $\E\in \Ext ^1(\I_\zeta\O (4K_X),\omega _X)$
sits in an exact sequence
$$0\raa\O (K_X+A)\ra\E\ra\O (K_X+D)\raa 0$$  
(since $\chi (\E)=\chi(K_X+A)+\chi (K_X+D)$).
But $\Ext ^1(\O (K_X+D),\O(K_X+A))=H^1(A-D)$.
Because $A+D\sim 3K_X$, we have $A-D\equiv K_X$. Now one can easily see
that $h^1(A-D)=0$ and $\E\simeq \O(K_X+A)\op \O(K_X+D)$.
Recall that  we have a surjection $\alpha  \:  \E\to \I_\zeta\O (4K_X)$.
Since 
$\O(-A)=\operatorname{im} \alpha |_{\O(K_X+D)}
\ot \O(-4K_X)\hookrightarrow \I_\zeta$,
we can assume that $A$ is an effective divisor and $\I_\zeta=\I_A+\I_D$. 

If we have two effective divisors
$D_1\equiv K_X$ and $D_2\equiv 2K_X$ intersecting at $\zeta$ and
$D_1+D_2\sim 3K_X$, then $\O(K_X+D_1)\op \O(K_X+D_2)$ gives  a non-trivial
extension of $\I_\zeta\O(4K_X)$  by $\omega_X$. Therefore from Serre's
construction $\zeta$ is contracted by $|4K_X|$, Q.E.D.
\medskip

\noindent {\bf 4.2.} {\it Proof of Theorem 0.2.}

From the preceding proposition $\zeta$ is contained in a curve 
$D\in |K_X+\tau|$, where $0\ne \tau\in \Tors X$.
One can easily see that $h^0(K_X+\tau)=1$ and $D$ is irreducible, reduced
(because $K_XD=1$) of genus $2$.  
Moreover, the trace of $|4K_X|$ on $D$ is a complete linear system.

If $2\tau=0$, then $\O_D(4K_X)=\O_D(2K_D)$ and because
the curve $D$ is honestly hyperelliptic we get 2 of the theorem. 

If $2\tau \ne 0$, then from the sequence
$$0\raa\O (K_X-2\tau)\ra\O(2K_X-\tau)\ra\O_D(2K_X-\tau)\raa 0$$
one gets $h^0(\O_D(2K_X-\tau))=1$ and $\zeta$ is a zero set of the unique
section of $\O_D(2K_X-\tau)$. Therefore  in $|2K_X-\tau|$ 
there is a unique divisor containing $D$ and 
all the other divisors in $|2K_X-\tau|$ intersect
$D$ exactly at  $\zeta$, Q.E.D.

\proclaim{Corollary 4.3}
Let $X$ be any canonical surface. Then $\varphi_{4K_X}$ is an embedding
if and only if $K_X^2\ge 2$ or $X$ is a torsion free Godeaux surface.
\endproclaim

The following theorem is a generalization of [Bo, Theorem 7.1]
to the case when the canonical model of a Godeaux surface
has singularities
(remark: the proof given in [Bo] cannot be easily generalized because
it uses the fact that if the image of $\varphi _{4K_X}(X)$ is singular
then $\varphi _{4K_X}$  is not biholomorphic at some points).

\proclaim{Corollary 4.4}
If $X$ is a Godeaux surface (i.e., a canonical numerical 
Godeaux surface with $\Tors X=\Z_5$), then $\varphi _{4K_X}$
contracts only $4$ tangent vectors (contained in $|K_X+\tau|$, 
$0\ne \tau \in \Tors X$) at the base points of $|2K_X|$. 
\endproclaim

\pr
By Theorem 0.2 $\varphi$ contracts only $4$ degree $2$ clusters,
which are scheme-theoretic intersections of curves
$D_1\in |K_X+\tau|$ and $2D_2$, where $D_2\in |K_X+2\tau|$.
But $D_1D_2=1$, hence $\zeta$ is a tangent vector at the point
$P=D_1\cap D_2$ and one can easily see that $P$ is a base point of $|2K_X|$, 
Q.E.D.

\heading 5. Numerical Campedelli surfaces \endheading

\proclaim{Theorem 5.1}
Let $C$ be a Gorenstein curve. If a degree $2$ cluster $\zeta$ 
is contracted by a linear system $|K_C|$ and each element of
$\Hom (I_\zeta,\O_C)$ is an injection, then
$C$ is honestly hyperelliptic.
\endproclaim

\pr
The proof is the same as the proof of [CFHR, Theorem 3.6].
\medskip

Let us also recall the following

\proclaim{Proposition 5.2} ([Re])
Let $X$ be a numerical Campedelli surface.
Then for every $\tau\in \Tors X$ we have $h^1(\tau)=0$. 
In particular, $h^0(K_X+\tau)=1$ for $\tau\ne 0$. 
\endproclaim

The proposition follows easily (by passing to the universal covering) 
from the fact that $q(X)=0$ and every \'etale Galois covering of
$X$ has bounded degree ($\le 9$, see [Re, Theorem I] or 
[Be, Remarque 5.9] with bound $\le 10$).

\medskip
\noindent {\bf 5.3.} {\it Proof of Theorem 0.3.}

The idea of the first part of the proof is stolen from [CFHR].

Let $C\in |\I_\zeta\O (2K_X)|$ and assume that $\zeta$ is contracted by
$|3K_X|$. Then we have a surjection
$$H^0(\O(3K_X))\ra H^0(\O_C(3K_X))=H^0(\omega _C),$$
so  $|\omega _C|$ also contracts  $\zeta$. Therefore
$\dim \Hom (\I_{\zeta,C},\O_C)=\dim \Hom (\I_\zeta\omega_C, \omega_C)=$

\noindent $h^1(\I_\zeta\omega_C)=h^1(\omega_C)+1=2$.

Because of Theorem 5.1 we can assume that there exists a nonzero
section $s\:\I_{\zeta,C}\to \O_C$, which is not an injection
(otherwise  $\varphi_{3K_X}|_C$ would be a double covering of $\P^1$). 
The section
$s$ vanishes on some subcurve $B\subset C$ and by the automatic adjunction
([CFHR, Lemma 2.4]) we have an injection $\I_\zeta\O_B(3K_X)\hookrightarrow
\omega _B$, which is generically a surjection. Therefore
$$\deg\I_\zeta\O_B(3K_X)=3K_XB-\deg (\zeta\cap B)\le \deg \omega_B,$$
i.e., 
$$3K_XB-\deg \omega_B\le \deg (\zeta\cap B)\le 2.$$
By the Bombieri connectedness theorem 
$C$ is not numerically 3-connected, $B\equiv K_X$ and $\deg (\zeta\cap B)=2$
(see the proof of [CFHR, Lemma 4.2]). But this means that $\zeta$
is contained in  $B$. Now it is sufficient to prove the following

\proclaim{Lemma 5.4}
In a notation as above $\zeta$ is a scheme-theoretic intersection of
two Cartier divisors  $B_1\in |K_X+\tau|$ and $B_2\in |K_X-\tau|$,
for some $\tau\in \Tors X$ such that $2\tau\ne 0$.
\endproclaim

\pr
There exists a reflexive sheaf  $\E$ sitting in the exact sequence
$$0\raa\omega_X\ra\E\mathop{\ra}^\alpha\I_\zeta\O(3K_X)\raa 0.$$
We already have one curve $B_1=B\equiv K_X$ containing $\zeta$. This gives 
an embedding $\O(3K_X-B_1)\hookrightarrow \I_\zeta\O(3K_X)$.
Because $\Ext ^1(\O(3K_X-B_1),\omega_X)=0$, this embedding lifts to
$\O(3K_X-B_1)\mathop{\hookrightarrow}^\beta \E$.
One can easily see that $\coker \beta$ is torsion-free
(otherwise there would exist a curve $B_1'<B$ containing
$\zeta$, which gives a contradiction). Moreover,
$(\coker \beta)^{**}=\O(K_X+B_1)$ and 
$\chi (\E)=\chi (K_X+B_1)+\chi (3K_X-B_1)$,
so $\coker \beta=\O(K_X+B_1)$. 
Now $\dim\Ext ^1(\O(K_X+B_1),\O(3K_X-B_1))=
h^1(\O(2\tau))=0$, so $\E\simeq \O(K_X+B_1)\op\O(3K_X-B_1)$.
The image of $\alpha|_{\O(K_X+B_1)}$ gives a divisor
$B_2=2K_X-B_1$ containing $\zeta$ and such that
$\I_{B_1}+\I_{B_2}=\I _\zeta$. If $B_1\in |K_X+\tau|$, then
$2\tau\ne 0$, because $h^0(K_X+\tau)=1$, Q.E.D.
\medskip

\noindent
{\sl Remarks.}

{(1)} Combining  Proposition 5.2 and Theorem 0.3 we see that there 
is only a finite
number of degree $2$ clusters which are contracted by $|3K_X|$ and 
are not contained in curves from 1 in Theorem 0.3.
It allows for a simple proof of the fact that $\varphi _{3K_X}$
is birational (pass to the 4-th point of the proof in [BC]).

{(2)} Theorem 0.3 could be proven by somewhat longer arguments but similar
as in the proof of  Theorem  0.2. Namely, if $\zeta$ is contracted by 
$|3K_X|$ one can construct a  rank  $2$ reflexive sheaf 
$\E\in  \Ext ^1(\I_\zeta \O(3K_X), \omega _X)$. A curve 
$C\in |\I_\zeta\O(2K_X)|$ defines an injection 
$\omega_X\hookrightarrow \I_\zeta \O(3K_X)$ lifting to  $\omega_X\to \E$,
since $\Ext ^1(\omega _X,\omega _X)=0$. Consider all the linear 
combinations $L$ of  our two  maps from $\omega_X$ to $\E$.  If a cokernel
of any of them has torsion then one can easily prove that we  are in  
case 2 of Theorem 0.3. Otherwise one can prove that the curve $C$   
is honestly hyperelliptic using the linear system  $L$.
In this  last case  the sheaf $\E$, which occurs to be a  bundle, 
is  stable and by Donaldson theorem from gauge theory it
(or rather the bundle $\E(-2K_X)$ which has trivial Chern classes)
corresponds to an irreducible $\hbox{SU} (2)$-representation  of $\pi_1(X)$
(in fact one should pull back $\E$ to the minimal model of $X$
since Donaldson's theorem holds for smooth surfaces; the reverse
is slightly harder).
Therefore  we have the following  proposition:

\proclaim{Proposition 5.5}
There is a bijection  between the set of honestly hyperelliptic
curves $C\in |2K_X|$ and the set of irreducible $\hbox{SU}(2)$-representations
of $\pi_1(X)$. 
\endproclaim

This proposition together with Theorem 0.3 shows that all the clusters
contracted  by $|3K_X|$  depend on the topology of $X$:
either they come from the torsion group $H_1(X, {\Z})$ or from
the representations of  the  fundamental  group.
This generalizes  [Ko, Theorem 2].

\medskip 
\noindent
{\bf 5.6.} {\sl Examples.} 

(1) Note that if $X$ is numerical Campedelli with 
$\pi_1^{alg}(X)=\Z_2\op\Z_2\op\Z_2$ (i.e., $X$ is a Campedelli surface)
then Proposition 5.5 together with Theorem 0.3 imply that 
$\varphi _{3K_X}$ is an embedding (since in this case 
$\pi_1(X)=\pi_1^{alg}(X)$ has no irreducible $\hbox{SU} (2)$-representations).
This was not known even though we knew an explicit description of $X$ 
(see [Pe, Remark after Theorem 2]).

(2) If $X$ is numerical Campedelli with 
$\pi_1^{alg}(X)=\H$ then there are no contracted clusters coming from 
torsion (see 1.6) but 
$\pi_1(X)=\H$ has an irreducible $\hbox {SU} (2)$-representation so
$\varphi _{3K_X}$ is not an embedding. 

\medskip
Summarizing all the known results, we get the following:

\proclaim{Corollary 5.7}
Let $X$ be any canonical surface. Then
$\varphi _{3K_X}$ is an embedding is an embedding if and only if
$K_X^2>2$ or $X$ is a Campedelli surface or
a numerical Campedelli surface with
$\pi_1^{alg} X=\Z_2\op\Z_2$, $\Z_2$ or $0$
and such that 
$\pi_1(X)$ has no irreducible $\hbox{SU} (2)$-representations.
\endproclaim

In view of this corollary  it would be very interesting to prove the following

\proclaim{Conjecture 5.8}
For any numerical Campedelli surface $\pi_1^{alg}(X)=\pi_1(X)$.
\endproclaim


\Refs
\widestnumber\key{CFHR}

\ref\key Be
\by A.~Beauville
\paper L'application canonique pour les surfaces de type general
\jour Invent. Math.
\vol 55 \yr 1979 \pages 121--140
\endref

\ref\key Bo
\by E. Bombieri
\paper The pluricanonical map of a complex surface
\paperinfo Lecture Notes in Math.
\vol 155 \publ Springer Verlag \yr 1970 \pages 35--87
\endref


\ref\key{BC} 
\by E.~Bombieri, F.~Catanese
\paper The tricanonical map of a surface with $K^2=2$, $p_g=0$
\paperinfo in C.~P.~Ramanujam, A Tribute, Tata Inst.
\publ Springer--Verlag  \yr 1978 \pages 279--290
\endref


\ref\key{CFHR} 
\by F. Catanese, M. Franciosi, K. Hulek, M. Reid
\paper Embeddings of curves and surfaces
\paperinfo preprint
\yr 1996 \endref

\ref\key{CT}
\by  F.~Catanese, F.~Tovena
\paper Vector bundles, linear systems and extensions of $\pi _1$
\paperinfo Lecture Notes in Math. \vol 1507 \publ Springer-Verlag \yr 1992 
\pages 51--71
\endref

\ref\key{Ha} 
\by R.~Hartshorne
\paper Stable reflexive sheaves
\jour Math.~Ann. \vol 254
\yr 1980 \pages 121--176
\endref

\ref\key{Ka} 
\by Y.~Kawamata
\paper Abundance theorem for minimal threefolds
\jour Invent.~Math. \vol 108 \yr 1992 \pages 229--246
\endref

\ref\key{Ko}
\by D.~Kotschick
\paper On the pluricanonical maps of Godeaux and Campedelli surfaces
\jour International J. of Math. \vol 5 \yr 1994 \pages 53--60
\endref

\ref\key La1 
\by A.~Langer
\paper Adjoint linear systems on normal surfaces
\paperinfo to appear in {J. Algebraic Geom}
\endref

 
\ref\key La2
\bysame
\paper Adjoint maps of algebraic surfaces
\paperinfo Ph.~D. Thesis
(in Polish), Warsaw University \yr 1998
\endref

\ref\key La3
\bysame
\paper A note on $k$-jet ampleness on surfaces
\paperinfo preprint \yr 1998
\endref


\ref\key Pe
\by C.~Peters
\paper On two types of surfaces of general type with
vanishing geometric genus
\jour Invent.~Math. \vol 32 \yr 1976 \pages 33--47
\endref

\ref\key{Re} 
\by M.~Reid
\paper Surfaces with $p_g=0$, $K^2=2$
\paperinfo unpublished manuscript and letters
\endref

\ref\key{We} 
\by L.~Weng
\paper A result on bicanonical maps of surfaces of 
general type
\jour Osaka J.~Math. 
\vol 32 \yr 1995 \pages 467--473
\endref


\ref\key{Xi} 
\by G.~Xiao
\paper Finitude de l'application bicanonique des surfaces
de type general
\jour Bull.~Soc. Math.~France
\vol 113 \yr 1985 \pages 23--51
\endref

\endRefs
\enddocument